\documentclass{amsart}

\usepackage{times}      		
\usepackage{latexsym}   	
\usepackage{amssymb}    	
\usepackage{amsmath}   	 
\usepackage{amsbsy}
\usepackage{amsgen}
\usepackage{amsfonts}
\usepackage{array}
\usepackage[all]{xy}    
\usepackage{epsfig}
\usepackage{hyperref}

\theoremstyle{plain}
\newtheorem{theorem}{Theorem}

\newtheorem*{thma}{Theorem A}
\newtheorem*{thmb}{Theorem B}
\newtheorem*{thmc}{Theorem C}

\newtheorem*{corob}{Corollary B} 
\newtheorem{corollary}[theorem]{Corollary} 
\newtheorem{lemma}[theorem]{Lemma} 
\newtheorem{proposition}[theorem]{Proposition}

\theoremstyle{remark}
\newtheorem*{ack}{Acknowledgements}
\newtheorem{remark}[theorem]{Remark}

\theoremstyle{definition}
\newtheorem{definition}[theorem]{Definition} 

\newtheorem{example}[theorem]{Example}

\def\lpq{L_{p, q}}

\def\rrr{\mathbb{R}}
\def\ccc{\mathbb{C}}
\def\zzz{\mathbb{Z}}

\DeclareMathOperator{\diam}{diam}

\def\bdm{\begin{displaymath}}
\def\edm{\end{displaymath}}
\def\beq{\begin{equation}}
\def\eeq{\end{equation}}
\def\bes{\begin{equation*}}
\def\ees{\end{equation*}}
\def\epcm{\end{picture}\end{center}\end{minipage}}
\def\bpcm{\begin{minipage}{80pt}\begin{center}\begin{picture}}

\def\t2{T^2}

\def\f4{F_4}
\def\g2{G_2}

\def\p2{\frac{\pi}{2}}
\def\dist{\textrm{dist}}

\def\txt{\textrm}

\def\Kl{\txt{Kl}}

\def\dim{\textrm{dim}}

\def\kl{Kl}

\begin{document}

\title[Cohomogeneity One Alexandrov Spaces]{Cohomogeneity One Alexandrov Spaces}

\title[Cohomogeneity One Alexandrov Spaces]{Cohomogeneity One Alexandrov Spaces}

\author[Galaz-Garcia]{Fernando Galaz-Garcia$^*$}
\address[Galaz-Garcia]{Instituto de Matem\'aticas, Universidad Nacional Aut\'onoma de M\'exico, Cuernavaca, Morelos, MEXICO}
\curraddr{Mathematisches Institut, WWU M\"unster, GERMANY}
\email{f.galaz-garcia@uni-muenster.de}

\author[Searle]{Catherine Searle$^{**}$}
\address[Searle]{Instituto de Matem\'aticas, Universidad Nacional Aut\'onoma de M\'exico, Cuernavaca, Morelos, MEXICO}
\email{csearle@matcuer.unam.mx}


\thanks {$^*$ The author was supported in part by CONACYT Project \#SEP-82471.}

\thanks{$^{**}$
The author was supported in part by 
CONACYT Project \#SEP-CO1-46274, CONACYT Project \#SEP-82471 and UNAM DGAPA project IN-115408.}

\date{\today}


\maketitle


\begin{abstract}
We obtain a structure theorem for closed, cohomogeneity one Alexandrov spaces and we classify closed, cohomogeneity one Alexandrov spaces in dimensions $3$ and $4$. As a corollary, we obtain the classification of
closed, $n$-dimensional, cohomogeneity one Alexandrov spaces admitting an isometric $T^{n-1}$ action.
 In contrast to the $1$- and $2$-dimen\-sional cases, where it is known that an Alexandrov space is a topological manifold, in dimension $3$ the classification contains, in addition to the known cohomogeneity one manifolds,  the spherical suspension of $\rrr P^2$, which is not a manifold.
\end{abstract}


\section*{Introduction}
The study of spaces with ``large" symmetry groups is of long-standing interest. One possible measure for the size of an isometric group action on a metric space is the dimension of the orbit space, called the \emph{cohomogeneity} of the action. From this point of view, a ``large'' action is an action of low cohomogeneity. Thus, a transitive action will be the largest possible symmetric action and 
 a cohomogeneity one action can then be considered as the second ``largest" symmetric action. Mostert introduced the concept of a cohomogeneity one action on a manifold in his 1956 paper ``On a compact Lie group acting on a Manifold" \cite{Mo}. Actions of cohomogeneity one are of interest in the theory of group actions on Riemannian manifolds and have been studied extensively in the presence of positive and non-negative sectional curvature, as well as in the curvature free setting (see, for example, \cite{Mo}, \cite{N}, \cite{P},  \cite{AA},  \cite{S},  \cite{V1}, \cite{V2}, \cite{GZ1}, \cite{GZ2},  \cite{GWZ}, \cite{H}).  In this article we propose to study isometric actions of cohomogeneity one on Alexandrov spaces, a class of metric spaces that properly contains the class of Riemannian manifolds. 

Alexandrov spaces play an important role in Riemannian geometry. They are a natural synthetic generalization of Riemannian manifolds with a lower curvature bound and the natural process of taking Gromov-Hausdorff limits is closed in Alexandrov geometry. 
The relationship between Riemannian and Alexandrov geometry has been used repeatedly to solve difficult problems in Riemannian geometry.
Perhaps the most well-known case is Perelman's solution to Thurston's Geometrization conjecture for $3$-manifolds, where Perelman's knowledge of a structure theorem for $3$-manifolds that collapse with a uniform lower curvature bound played a crucial role (cf. \cite{Pe1}, \cite{Pe2}, \cite{Pe3}, \cite{SY}).
As in Riemannian geometry, it is natural that in order to try to understand Alexandrov spaces in general one begins by understanding those with lots of symmetry.

A well-known result of Myers and Steenrod \cite{MS} states that the group of isometries of a Riemannian manifold is a Lie group.
In 1994 Fukaya and Yamaguchi  \cite{FY}  proved that the isometry group of
a length space which is locally compact and of finite Hausdorff dimension, with curvature bounded away from $-\infty$, will be a Lie group. It follows that the isometry group of an Alexandrov space is a Lie group. Moreover, if we restrict our attention to isometric actions on closed Alexandrov spaces,
we need only consider compact Lie groups. As a consequence,  a homogeneous closed Alexandrov space is in fact a Riemannian manifold. It is then a natural next step to study closed Alexandrov spaces of cohomogeneity one.

In the case of a cohomogeneity one action on a closed Alexandrov space, the only possible orbit spaces are a closed interval or a circle. We extend  structural results about cohomogeneity one manifolds to the Alexandrov setting and classify cohomogeneity one, closed Alexandrov spaces in dimensions $3$ and $4$, extending the classification results of Mostert \cite{Mo} and Neumann \cite{N}, in dimension $3$, and of Parker \cite{P}, in dimension 4, noting that in his thesis,  Hoelscher \cite{H}
corrects an omission in Parker's classification.

Many results about cohomogeneity one manifolds are of a group theoretic nature. This is the case, for example, for theorems about  the basic decomposition of a cohomogeneity one manifold with quotient space an interval as a union of disc bundles, its fundamental group, reductions and extensions, and primitive and non-primitive actions. We discuss their generalization to Alexandrov spaces in Section \ref{S:3}.

We remark that    we cannot always extend results related to the smooth structure of a cohomogeneity one manifold, such as those associated with the \emph{Weyl group} of a  cohomogeneity one Riemannian $G$-manifold $M$. The Weyl group, $W$, of $M$ is defined as the stabilizer of a normal geodesic, $C$, perpendicular to all orbits of the action, modulo its kernel $H$ (see \cite{GWZ}). It is known (cf. \cite{AA}) that $C$ is either a circle or a one-to-one immersed line
and $M/G=C/W$, where $W$ is a dihedral subgroup of $N(H)/H$ and $N(H)$ is the normalizer of $H$. Further, when $W$ is finite, $C$ will be closed, that is, a circle. In the Alexandrov case, however,
we may have $N(H)/H=\{1\}$ with $C$ a closed geodesic segment. This is the case, for example, when we consider the suspension of $\rrr P^2$  with the suspension of the standard transitive $SO(3)$ action on $\rrr P^2$. In this case a normal geodesic perpendicular to all the orbits terminates in the singular orbits, where the space of directions is not a sphere.

 Our main theorem is the  following structure result for closed Alexandrov spaces of cohomogeneity one, which generalizes the structure of a closed, cohomogeneity one manifold. We denote the cone over a topological space $Y$  by $C(Y)$, and the spherical suspension of $Y$ by $\Sigma(Y)$.


\begin{thma} Let $X$ be a closed Alexandrov space with an effective isometric $G$ action of cohomogeneity one with principal isotropy $H$. 
 \begin{itemize}
 	\item[(1)]  If the orbit space of the action is an interval, then $X$ is the union of two fiber bundles over the two singular orbits whose fibers are cones over positively curved homogeneous spaces, that is,
	\[ 
		X=G\times_{K_{-}}C(K_{-}/H) \cup_{G/H}G\times_{K_{+}} C(K_{+}/H).
	\]
	The group diagram of the action is given by $(G, H, K_{-}, K_{+})$, where $K_{\pm}/H$ are positively curved homogeneous spaces. Conversely, a group diagram $(G, H, K_{-},K_{+})$, where $K_{\pm}/H$ are positively curved homogeneous spaces, determines a cohomogeneity one Alexandrov space.
	\\
  	\item[(2)] If the orbit space of the action is a circle, then $X$ is equivariantly homeomorphic to a $G/H$-bundle over a circle with structure group $N(H)/H$. In particular, $X$ is a manifold.
	\end{itemize}
	
\end{thma}

As mentioned earlier, a homogeneous Alexandrov space is a homogeneous manifold. In view of this, it is natural to ask at what level of symmetry can Alexandrov spaces occur which are no longer manifolds. Since $1$- and $2$-dimensional Alexandrov spaces are topological manifolds \cite{BBI}, it follows  from Mostert \cite{Mo} that a closed $1$- or $2$-dimensional Alexandrov space of any cohomogeneity is equivariantly homeomorphic to a closed manifold with the same action.  In contrast to the $1$- and $2$-dimensional cases, a cohomogeneity one $3$-dimensional Alexandrov space will not necessarily be a manifold, as the following theorem shows. 


\begin{thmb} Let $X^3$ be a closed, $3$-dimensional Alexandrov space with an effective isometric cohomogeneity one $G$ action. Then $G$ is $SO(3)$ or $T^2$ and the only non-manifold we obtain is $\Sigma(\rrr P^2)$.

\end{thmb}

Theorem B can be obtained from Theorem~A by analyzing all possible admissible group diagrams, or by directly analyzing the 
spaces of directions at the singular orbits.

The following corollary is an immediate consequence of Theorem~B. We will denote by $S^2\tilde{\times}S^1$ the non-trivial $S^2$ bundle over $S^1$, by $\lpq$ a lens space, by $\Kl$ a Klein bottle, by $Mb$ the compact M\"obius band and by $A$ the manifold $Mb \times S^1 \cup S^1\times Mb$ intersecting canonically in $S^1\times S^1$. 


\begin{corob} A closed Alexandrov space $X^n$ of dimension $n\geq 4$ with an effective isometric $T^{n-1}$ action is equivariantly homeomorphic to the product of  $T^{n-3}$ with one of $T^3$, $S^3$, $\lpq$, $S^2\times S^1$, $S^2\tilde{\times} S^1$, $\kl\times S^1$, $\rrr P^2\times S^1$ or   $A$.
\end{corob}

\begin{remark} If we allow the orbit space to be non-compact, we obtain a complete generalization of Mostert's original result for compact Lie groups. In particular, it follows  from the proof of Theorem B that the only $3$-di\-men\-sional Alexandrov space of cohomogeneity one with non-compact orbit space that is not a manifold is the cone over $\rrr P^2$.
\end{remark}
\smallskip

In dimension $4$ the number of cohomogeneity one Alexandrov spaces that are not manifolds increases substantially. 
Including the suspensions of the $3$-dimensional homogeneous spherical space forms, we find $15$ classes of group diagrams corresponding to cohomogeneity one Alexandrov spaces  that are not manifolds. 
Our last main result, stated in Section~\ref{S:D4}, is the classification of closed, cohomogeneity one $4$-dimensional Alexandrov spaces, extending Parker's classification of closed, cohomogeneity one $4$-dimensional manifolds (cf. \cite{P}, \cite{H}).

Finally, we note that simply connected cohomogeneity one manifolds in dimensions up to and including $7$ have been classified (\cite{Mo}, \cite{N}, \cite{P}, \cite{H}), but the non-simply connected ones have only been classified through dimension $4$. In order to continue the classification of cohomogeneity one Alexandrov spaces in higher dimensions, one must first classify the non-simply connected cohomogeneity one manifolds in dimensions $5$ and above. Moreover,  a classification of non-simply connected homogeneous spaces of positive curvature for dimensions greater than $5$ is notably lacking in the literature, and this knowledge is also required to complete the classification of cohomogeneity one Alexandrov spaces in higher dimensions.

The present article is organized as follows. In Section~\ref{S:2} we recall some basic results from Alexandrov geometry and about cohomogeneity one actions and prove Theorem~A. In Section~\ref{S:3} we extend several group theoretic results on cohomogeneity one manifolds to the Alexandrov setting. Section~\ref{Proof_thm_B} contains the proof of Theorem~B and the proof of Corollary~B is presented in Section~\ref{Proof_Cor_B}. Closed, cohomogeneity one Alexandrov spaces of dimension $4$ are classified in Section~\ref{S:D4}.


\begin{ack} This work was done during the first-named author's tenure as a Postdoctoral Fellow at the Institute of Mathematics of the UNAM in Cuernavaca, Mexico. Both authors would like to thank Karsten Grove for helpful conversations and for suggesting improvements to Theorem A. The authors also would like to thank one of the referees for suggesting improvements to Proposition~\ref{p:normal_space}.
\end{ack}

\section{Preliminaries and Proof of Theorem~A}
\label{S:2}


We first fix some notation. We will denote an Alexandrov space by $X$. Given an isometric (left) action $G\times X\rightarrow X$ of a Lie group $G$, and a point $x\in X$, we let $G(x)=\{\,gx :g\in G \,\}$ be the \emph{orbit}  of $x$ under the action of $G$. The \emph{isotropy group} of $x$ is the subgroup $G_x=\{\, g\in G: gx=x\,\}$. Observe that $G(x)\simeq G/G_x$. We will denote the orbit space of this action by $X^*=X/G$. Similarly, the image of a point $x\in X$ under the orbit projection map $\pi:X\rightarrow X^*$ will be denoted by $x^*\in X^*$. We will assume throughout that $G$ is compact and its action is \emph{effective}, i.e., that $\bigcap_{x\in X}G_x$ is the trivial subgroup $\{e\}$ of $G$, except in Section~\ref{S:D4}, where we will assume that the action is \emph{almost effective}, i.e., that  $\bigcap_{x\in X}G_x$ is a finite subgroup of $G$.


\subsection{Alexandrov geometry} We recall some basic facts about Alexandrov spaces, all of which can be found in \cite{BBI}. A finite dimensional length space $(X,\mathrm{dist})$  has curvature bounded from below by $k$ if every point $x\in X$ has a neighborhood $U$ such that for any collection of four different points $(x_0,x_1,x_2,x_3)$ in $U$, the following condition  holds:
\[
\angle_{x_{1},x_{2}}(k)+\angle_{x_{2},x_{3}}(k)+\angle_{x_3,x_1}(k)\leq 2\pi.
\]
Here $\angle_{x_{i},x_{j}}(k)$, called the \emph{comparison angle}, is the angle at $x_{0}(k)$ in the geodesic triangle in $M^2_k$, the simply-connected $2$-manifold with constant curvature $k$, 
with vertices $(x_{0}(k),x_{i}(k),x_{j}(k))$, which are the isometric images of $(x_{0},x_{i},x_{j})$. An Alexandrov space is a complete length space of curvature bounded below by $k$, for some $k\in \rrr$.

The \emph{space of directions} of a general Alexandrov space $X^n$ of dimension $n$ at a point $x$ is,
by definition, the completion of the 
space of geodesic directions at $x$. We will denote it by $\Sigma_xX^n$.
It is a compact Alexandrov space of dimension $n-1$ with curvature bounded below by 1. 

We recall the following result (cf. \cite{GW}, Lemma~2.5) which will be used in the proof of Proposition \ref{p:normal_space}.

\begin{lemma}[Join Lemma]\label{l:join} Let $X$ be an $n$-dimensional Alexandrov space with curvature bounded below by $1$. If $X$ contains the unit round sphere $S^m_1$ isometrically, then 
$E=\{\,x\in X : \mathrm{dist}(x,S_1^m)=\pi/2\,\}$ is an isometrically embedded $(n-m-1)$-dimensional Alexandrov space with curvature bounded below by $1$, and $X$ is isometric to $S_1^m*E$ with the standard join metric  (cf. \cite{GP,GM}).
\end{lemma}



\subsection{Cohomogeneity one Alexandrov spaces}
It is well-known that the orbit space of a closed manifold of cohomogeneity one is either a closed interval or a circle. 
A cohomogeneity one $G$ action on a closed manifold with orbit space an interval determines a group diagram 

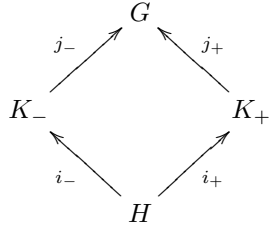
\begin{equation*}
\xymatrix{ & G    & \\
K_{-} \ar[ru]^{j_{-}} & & K_{+} \ar[lu]_{j_{+}}\\
& H  \ar[lu]^{i_{-}} \ar[ru]_{i_{+}}& }
\end{equation*}
 where $i_{\pm}$ and $j_{\pm}$ are the inclusion maps, $K_{\pm}$ are the isotropy groups of the singular orbits at the endpoints of the interval, and $H$ is the principal isotropy group of the action. 
 We will denote this diagram by the $4$-tuple $(G, H, K_{-}, K_{+})$. 
 We remark that the inclusion maps are an important element in the group diagram and together with the group diagram determine the cohomogeneity one space. The following example illustrates this point: $(T^2, \{e\}, T^1, T^1)$ determines both $S^3$ and $S^2\times S^1$, where in the first case the inclusion maps are to the first and second factors, respectively, of $T^2$, and in the second case, both inclusion maps are the same (cf. \cite{N}).

 It is well-known that $K_{\pm}/H$ are spheres and that $M$ decomposes as a union of disc bundles over the corresponding singular orbits, $G/K_{\pm}$. 
 Further, there is a one-to-one correspondence between  group diagrams $(G,H,K_{-},K_{+})$ with $K_{\pm}/H$ a sphere, and cohomogeneity one manifolds (cf. \cite{Mo}). 
In the case where the orbit space is a circle, $M$ is a  $G/H$-fiber bundle over $S^1$ with structure group $N(H)/H$. 

 Parker \cite{P} showed that in the special case where $K_{+}\supseteq K_{-}\supset H$, 
a cohomogeneity one manifold $M$ given by $(G, H, K_{-}, K_{+})$ is a fiber bundle over $G/K_{+}$ with fiber a cohomogeneity one
submanifold with group diagram $(K_{+}, H, K_{-}, K_{+})$.
This is part of a more general phenomenon called \textit{non-primitivity}. We recall that a cohomogeneity one action is called \textit{non-primitive} if there is a $G$-equivariant map
$X\rightarrow G/L$ for some subgroup $L\subset G$ (see \cite{AA}, p.17). Otherwise, the 
action is said to be \emph{primitive}. For cohomogeneity one manifolds, non-primitivity is equivalent 
to the statement that, for some group diagram representation $(G,H,K_{-},K_{+})$, we have $H \subset \{K_{-}, K_{+}\}\subset L\subset G$, i.e., for some 
invariant metric and some normal geodesic, $K_{\pm}$ generate a proper subgroup of $G$. In terms of
the original groups, the action is primitive if $K_{-}$ and $nK_{+}n^{-1}$ generate 
$G$, for any fixed $n\in  N(H)_c$, where $c$ is a normal geodesic in $M$. 

In the case  of a cohomogeneity one Alexandrov space we obtain analogous structure results. As in the manifold case, the orbit space of a closed cohomogeneity one Alexandrov space is either a circle or an interval, and a closed Alexandrov space with orbit space an interval decomposes as the union of two bundles over the singular orbits. In contrast to the case  of a cohomogeneity one manifold, the Alexandrov case is less rigid. The bundles over the singular orbits need not be disc bundles, but will have fibers that are cones over positively curved homogeneous spaces (which are discs when the homogeneous space is a sphere). There is also a one-to-one correspondence between admissible group diagrams and cohomogeneity one Alexandrov spaces.
 In the special case where $K_{+}\supseteq K_{-}\supset H$, a cohomogeneity one Alexandrov space will be a bundle over $G/K_{+}$ with fiber a cohomogeneity one Alexandrov subspace with group diagram
$(K_{+}, H, K_{-}, K_{+})$.
In the case of equality the fiber of $(K_{+}, H, K_{+}, K_{+})$ is $\Sigma(K_{+}/H)$. In the particular case of an Alexandrov space, we can only talk about normal geodesic segments, but this is sufficient to then be able to say that  in terms of the original groups, the action is primitive if $K_{-}$ and $nK_{+}n^{-1}$ generate 
$G$, for any fixed 
$n\in  N(H)_c$, where $c$ is a normal geodesic segment in $X$.


\subsection{Proof of Theorem A} We prove assertions (1) and (2) separately as Propositions~\ref{p:interval} and \ref{p:circle}, respectively. We begin by defining the set of normal directions to a subset of the space of directions.


\begin{definition} Let $X$ be an Alexandrov space and fix $x\in X$. Given $A\subset \Sigma_xX$, we define the set of \emph{normal directions} to $A$, denoted by $\nu(A)$, by
\[
\nu(A)=\{\, v\in \Sigma_xX : \dist(v,w)=\diam(\Sigma_xX)/2\text{ for all } w\in A\,\}.
\]
\end{definition}

Note that when $X$ is a Riemannian manifold, this definition coincides with the usual definition of normal directions.
In the specific case of an Alexandrov  space $X$ admitting an isometric $G$ action, we can express the space of directions  at any point $x\in X$ as a join of the unit tangent space to the orbit $G(x)$, denoted by $S_x$, and its normal  space $\nu(S_x)$, provided $\dim(G/G_x)>0$, as we see in the following proposition.





\begin{proposition}\label{p:normal_space} Let $X$ be an Alexandrov space admitting an isometric $G$ action and fix $x\in X$ with $\dim(G/G_x)>0$. If $S_x\subset \Sigma_xX$ is the unit tangent space to the orbit $G(x)\simeq G/G_x$,  then the following hold. 
\\
\begin{itemize}
	\item[(1)] The set $\nu(S_x)$ is a compact, totally geodesic Alexandrov subspace of $\Sigma_xX$ with curvature bounded below by $1$,  and the space of directions $\Sigma_x X$ is isometric to the join $S_x * \nu(S_x)$ with the standard join metric.\\
	\item[(2)] Either $\nu(S_x)$ is connected or it contains exactly two points at distance $\pi$.\\
\end{itemize} 
\end{proposition}


\proof 
Compactness of $\nu(S_x)$ follows from the continuity of the distance function. The rest of  part (1) will follow  from the Join Lemma for Alexandrov spaces (
see 
Lemma \ref{l:join}), once we show that  $S_x$ is isometric to the 
round unit sphere and
 totally geodesic in $\Sigma_xX$. 
The fact that the unit tangent sphere $S_x$ to the orbit $G(x)$ is isometric to the unit round sphere is a consequence of  the tangent cone  at $x$ isometrically splitting as $\mathbb{R}^n \times U$, where $\mathbb{R}^n$, 
the tangent space to the orbit $G(x)\simeq G/G_x$
, is equipped with the standard Euclidean metric, and $U$ is a cone.


Let $v$ be a basis element of the tangent space of $G/G_x$. To see that the tangent cone $C_xX$ contains a line passing through $-v$, the cone point, and $v$, consider a $1$-parameter subgroup $c(t) = \exp(vt)x$ of the orbit $G/G_x$.
Take a sequence $t_n \longrightarrow 0$ and,  for each $n$,  rescale the metric on $X$ so that the distance between $c(0)$ and $c(t_n)$ is $1$. Let  $X_n$ denote $X$ with this rescaled metric. For each $n$, $X_n$ is equipped with an isometry $i_n$ sending $c(0)$ to $c(t_n)$. The sequence $X_n$ converges to the tangent cone $C_xX$ in the Gromov-Hausdorff sense, and the sequence of isometries $\{i_n\}$ converges to an isometry $i:C_x X \rightarrow C_x X$ sending the cone point $o_x$ to $v$, the direction of $c(t)$, which implies there is a line. Indeed, this orbit is a line and not, 
for example, 
a circle (which is a possibility for the orbit of a point of a one parameter family of isometries on $\mathbb{R}^n$), since it is constructed as a blow up limit. Observe that  the line constructed is also the line 
passing 
through $-v$, the cone point, and $v$. By the Splitting Theorem~\cite{M}, the cone must split isometrically as $\mathbb{R} \times U$. Doing this for all basis elements of the tangent space of $G/G_x$  splits off  $\mathbb{R}^n$.

To prove part (2), observe first that $\diam\Sigma_xX=\pi$. Arguing as in \cite{GG} one can see that $\nu(S_x)$ is \emph{totally $\pi$-convex}, i.e., a geodesic of length less than $\pi$ between any two points in $\nu(S_x)$ is entirely contained in $\nu(S_x)$. From this it follows that any two connected components of $\nu(S_x)$ are at distance at least $\pi$ apart. Standard triangle comparison arguments then show that there can be at most two such components and in that event they must be points.

\qed




\begin{proposition}\label{p:interval} If $X$ is a closed, cohomogeneity one Alexandrov space with orbit space an interval, then $X$ is the union of two bundles whose fibers are cones over positively curved homogeneous spaces, that is
 
	\[
	G\times_{K_{-}}C(K_{-}/H) \cup_{G/H}G\times_{K_{+}} C(K_{+}/H),
	\]
	where the group diagram of the action is given by $(G, H, K_{-}, K_{+})$ and $K_{\pm}/H$ are positively curved homogeneous spaces.
	Conversely, any diagram $(G,H,K_{-},K_{+})$, with $K_{\pm}/H$ positively curved homogeneous spaces, gives rise to a cohomogeneity one Alexandrov space. 
\end{proposition}
 \proof
 In this case there are two singular orbits with isotropy $K_\pm$. These orbits are manifolds $G/K_{\pm}$. Fix $x\in G/K_{\pm}$ so that $x^*$ is one of the endpoints of the orbit space $X^*$. 
 The tangent cone $K_xX$ contains the tangent space  $T_x\,G(x)$. Hence  the space of directions $\Sigma_xX$ contains the unit tangent sphere $S_x\subset T_xG(x)$. It follows from Proposition~\ref{p:normal_space} that $\nu(S_x)$ is a compact Alexandrov space with curvature $\geq 1$. Moreover, the isotropy group $K_{\pm}$ acts by isometries on $\Sigma_xX$, so it acts on $\nu(S_xX)$ and, since the space of directions of $X^*$ at an endpoint $x^*$ is a  point, it follows that the $K_{\pm}$ action on  $\nu(S_xX)$ is transitive and thus $\nu(S_xX)$ is a homogeneous space. 
 The desired bundle decomposition is then a consequence of  the Slice Theorem (cf. \cite{Br}). 
We note that in the case of an exceptional or a codimension one singular orbit, the cone bundle over the corresponding orbit is always a disc bundle as in the Riemannian case, since the only possibilities for the fiber 
in these two cases are $S^0$ and $S^1$, respectively.

 As shown above, if we have an Alexandrov space of cohomogeneity one we will obtain a group diagram $(G, H, K_{-}, K_{+})$, where $K_{\pm}/H$ are positively curved homogeneous spaces. To prove the converse, suppose we have group inclusions $H\leq K_{\pm} \leq G$ so that $K_{\pm}/H$ admit  homogeneous metrics with positive sectional curvature. We  fix such metrics with curvatures greater than or equal to $1$, and consider the corresponding spherical cones, $CK_{\pm}/ H$, which have  curvatures bounded below by $1$. Now, equip $G$ with the bi-invariant metric, so that the products $G \times CK_{\pm}/H$ are Alexandrov spaces with curvature bounded below by $0$.

The cone bundles,  $(G \times CK_{-}/H)/K_{-}$ and $(G\times  CK_{+}/H)/K_{+}$, are Alexandrov spaces with boundary $G/H$, but with possibly different metrics, $g_{-}$ and $g_{+}$, respectively, on their boundaries.
We consider the 
manifold, $G/H \times [0,1]$, with 
any $G$-invariant metric which is a product metric near its two boundary components, and with metrics $g_{-}$, and $g_{+}$, respectively, on these boundary components.
Observe that $G/H \times [0,1]$ is a Riemannian manifold with totally geodesic boundary components, so it is also an  Alexandrov space. By the Gluing Theorem (cf. \cite{Pet}) we can glue the three Alexandrov spaces, $(G \times CK_{-}/H)/K_{-}$, $(G\times  CK_{+}/H)/K_{+}$ and $G/H \times [0,1]$, together to form an Alexandrov space with an isometric $G$ action of cohomogeneity  one.
Observe that  in general the curvatures in $G/H \times [0,1]$ may change sign.
\qed


  
\begin{remark} Positively curved, simply connected homogeneous spaces  have been classified in all dimensions by Berger \cite{B}, Wallach  \cite{Wa},  Alloff and Wallach \cite{AW},  and Berard-Bergery \cite{BB}. There is a notable lack in the literature, however, of a classification of positively curved homogeneous spaces regardless of fundamental group.
Clearly, the non-simply connected ones are covered by those that have been classified. Combining this with the classification of homogeneous space forms due to Wolf \cite{W}, and the fact that in even dimensions there can at most be $\zzz_2$ quotients,  by Synge's Theorem, it follows that the positively curved homogeneous spaces in dimensions $5$ and below are $S^0$, $S^1$, $S^2$, $\rrr P^2$, the $3$-dimensional spherical space forms, $S^4$, $\rrr P^4$, $\ccc P^2$ (noting that $\ccc P^2$ admits no $\zzz_2$ quotient) and, finally, the $5$-dimensional spherical space forms.

 We can then  describe the space of directions  $\Sigma_{x_{\pm}}X^k$ for a $k$-dimensional cohomogeneity one $G$-Alexandrov space $X^k$ at a singular orbit $Gx_{\pm}$, for $2\leq k\leq 6$. We recall first that the homogeneous spherical space forms in dimension $k=3$ and $k=5$ are given by $S^k/\Gamma$, where $\Gamma$ is one of 
$\zzz_n$, $n\geq 1$, $D_{m}^*$, $m\geq 2$, $T^*$, $I^*$, $O^*$; here $^*$ denotes the corresponding binary group, and $T$, $I$, $O$ are, respectively, the tetrahedral, icosahedral and octahedral groups (see \cite{W} p. 89). 
 It follows that, in dimension $3$, $\Sigma_{x_{\pm}}X^3$ is one of $S^2$ and $\rrr P^2$; in dimension $4$, $\Sigma_{x_{\pm}}X^4$ is one of $S^3$, $SU(2)/\Gamma$  or $\Sigma(\rrr P^2)$; in dimension $5$, $\Sigma_{x_{\pm}}X^5$ is one of $S^4$, $\rrr P^4$, $\ccc P^2$, $\Sigma(S^3/\Gamma)$ and $\Sigma(\Sigma(\rrr P^2))=S^1*\rrr P^2$; finally, in dimension $6$, the possibilities for $\Sigma_{x_{\pm}}X^6$ are $S^5$, $\rrr P^5$,  $S^5/\Gamma$, $\Sigma(\rrr P^4)$, $\Sigma(\ccc P^2)$, $\Sigma(\Sigma(SU(2)/\Gamma))=S^1*SU(2)/\Gamma$ and $\Sigma(\Sigma(\Sigma(\rrr P^2)))=S^2*\rrr P^2$.
\end{remark}
  
Proceeding as in \cite{P}, we  obtain the following corollary.
  
  
 \begin{corollary} \label{c:interval} Let $X$ be a closed, cohomogeneity one Alexandrov space with orbit space an interval. In the special case where 
 $K_{+}\supseteq K_{-}\supset H$, 
 $X$ will decompose as 
 a fiber bundle over  $G/K_{+}$ with fiber the cohomogeneity one Alexandrov space with orbit space an interval and group diagram 
 $(K_{+}, H, K_{-}, K_{+})$.
 \end{corollary}

 We conclude this section with the case where the orbit space of $X$ is a circle.
 
  \begin{proposition}\label{p:circle} Let $G$ act on $X$, a closed, Alexandrov space, by  cohomogeneity one with orbit space a circle. Then $X$ is equivariantly homeomorphic to a fiber bundle over a circle with fiber a homogeneous space $G/H$ and structure group $N(H)/H$. 
  \end{proposition}
 
 \proof This follows immediately once we show that all orbits must be principal and thus homogeneous spaces. Suppose instead that we have a singular orbit, $G/K$. In this case, the singular isotropy group $K$ will act on the normal space of directions to the orbit with isotropy group $H$. The dimension of the normal space of directions is equal to the dimension of $K/H$ and thus the action must be transitive. In particular, this tells us that in the orbit space $X^*$, the normal space of directions is a point and thus the image of $G/K$ in $X^*$ is an endpoint. This yields a contradiction and thus all orbits are principal.  
\qed
  
 
\section{Results related to the group diagram}
\label{S:3}
A number of results about cohomogeneity one manifolds related to the group diagram of the action readily generalize to the Alexandrov space setting. We present them in this section.

Recall that a group diagram $(G,H,K_{-},K_{+})$, with $K_{\pm}/H$ positively cur\-ved homogeneous spaces, is in one-to-one correspondence with 
a closed, cohomogeneity one Alexandrov space. It is important to know when two different diagrams
can determine the same Alexandrov space. As in the manifold case, we say that the action of $G_1$ on $X_1$ is \textit{equivalent} to the action of $G_2$ on $X_2$ if there exists a homeomorphism $f:X_1\rightarrow X_2$ and an isomorphism $\phi:G_1\rightarrow G_2$ such that $f(g\cdot x)=\phi(x)\cdot f(x)$ for all $x\in X_1$ and $g\in G_1$.
In the special case where $G_1=G_2$ we can ask that the corresponding manifolds be $G$-equivariantly homeomorphic. In this case, the following proposition, obtained as in \cite{GWZ}, determines when two different group diagrams yield the same Alexandrov space.

\begin{proposition} If a cohomogeneity  one Alexandrov is given by a group diagram $(G,H,K_{-},K_{+})$, then any of the following operations on the group diagram will result in a $G$-equivariantly homeomorphic Alexandrov space:
\begin{enumerate}
 	\item Switching $K_{-}$ and $K_{+}$,
	\item Conjugating each group in the diagram by the same element of $G$,
	\item Replacing $K_{-}$ with $aK_{-}a^{-1}$  for $a\in N(H)_0$.
\end{enumerate}
Conversely, the group diagrams for two $G$-equivariantly homeomorphic cohomogeneity one, closed Alexandrov spaces must be mapped to each other by some combination of these three operations.
\end{proposition}

The following application of van Kampen's Theorem  is obtained as in \cite{H}.

\begin{proposition} Let $X$ be the cohomogeneity one Alexandrov space given by the group diagram $(G,H,K_{-},K_{+})$ with $K_{\pm}/H$ a positively curved, simply connected homogeneous space. Then $\pi_1(X)\cong \pi_1(G/H)/N_{-}N_{+}$, where
\[
N_{\pm}=\ker\{\,\pi_1(G/H)\rightarrow \pi_1(G/K_{\pm})\,\}=\mathrm{Im}\{\,\pi_1(K_{\pm}/H)\rightarrow \pi_1(G/H)\,\}.
\]
In particular $X$ is simply connected if and only if the images of $K_{\pm}/H$ generate $\pi_1(G/H)$ under the natural inclusions.
\end{proposition}

 The following three results about reductions and extensions of cohomogeneity one actions are also obtained as in \cite{H}.

  \begin{proposition}\label{p:red} Let $X$ be the cohomogeneity one Alexandrov space given by the group diagram 
$(G, H, K_{-}, K_{+})$ and 
suppose $G = G_1 \times G_2$ with $\pi_2 (H ) = G_2$, where $\pi_2$ is the natural projection onto $G_2$. Then the subaction of 
$G_1\times  1$ on $X$ is also by cohomogeneity one, with the same orbits, and with isotropy groups 
$K_{\pm}=K_{\pm}\cap (G_1\times 1)$ and $H_1=H\cap (G_1\times q)$.
\end{proposition}

For the next proposition we will need the concept of a \emph{normal extension}, which we now recall. Let $X$ be a cohomogeneity one Alexandrov space with group diagram $(G_1, H_{1},K_{-}, K_{+})$ and let $L$ be a compact, connected subgroup of $N(H_1)\cap N(K_{-})\cap N(K_{+})$. Observe that the subgroup $L\cap H_1$ is normal in $L$ and define $G_2:=L/(L\cap H_1)$. We can then define an action on $X$ by $G_1\times G_2$ orbitwise by $(\hat{g_1}, [l])\star g_1(G_1)_x=\hat{g_1}g_1l^{-1}(G_1)_x$
on each orbit $G_1 / (G_1)_x$ , for $(G_1)_x = H_1$ or $K_{\pm}$.
Such an extension is called a \textit{normal extension} of $G_1$.

\begin{proposition}\label{p:ext} A normal extension of $G_1$ describes a cohomogeneity one action of $G := G_1\times G_2$ on $X$ with the same orbits as $G_1$ and with group diagram 
$(G_1\times G_2, (H_1\times 1)\cdot \Delta L, (K_{-}\times 1) \cdot \Delta L, (K_{+}\times 1)\cdot \Delta L)$, 
where $\Delta L=\{(l, [l]): l\in L\}$.
\end{proposition}
  
  As in the case for manifolds, extensions and reductions are ``inverses" of each other, as the following proposition shows.
  
  \begin{proposition}\label{p:red-ext} 
 For $X$ as in Proposition \ref{p:red}, the action by $G = G_1\times G_2$ occurs as the normal extension of the reduced action of $G_1\times 1$ on $X$.
 \end{proposition}

\section{Proof of Theorem B} 
\label{Proof_thm_B}

We begin by recalling that for a compact Riemannian manifold the connected components of the fixed point set of an isometric circle action will be closed, totally geodesic submanifolds (cf. \cite{Ko}).  
For a $3$-dimensional Alexandrov space we obtain the following proposition.


\begin{proposition}\label{r:1} 
The connected components of the fixed point set of an isometric $S^1$ action on a closed $3$-dimensional Alexandrov space $X^3$ are intervals or circles.
\end{proposition}

\proof
Let $F$ be the fixed point set of the $S^1$ action and observe that $F$ must be a compact subset of $X^3$. We will show that the connected components of $F$ are topological $1$-manifolds, possibly with boundary. 

Given $x\in F$, the space of directions $\Sigma_xX^3$ is homeomorphic the sphere $S^2$ or to $\rrr P^2$ and the $S^1$ action on $\Sigma_xX^3$ fixes two points if $\Sigma_xX^3$ is $S^2$, or a single point if $\Sigma_xX^3$ is $\rrr P^2$. Observe now that a fixed direction in $\Sigma_xX^3$ corresponds to a half-line in the tangent cone $K_xX^3$ and two fixed directions correspond to a line in $K_xX^3$.  
The desired conclusion then follows from the existence of a homeomorphism between the tangent cone at $x$ and a neighborhood of $x$ in $X^3$ (cf. \cite{Per, K}).
\qed


It is well-known that the fixed point set components of a circle action on a compact smooth $3$-manifold are circles. The following example illustrates  Proposition~\ref{r:1} when the fixed point set component is a closed interval. 

 \begin{example} Recall that $\rrr P^2$ equipped with its standard positively curved Riemannian metric admits an isometric $S^1$ action with a single fixed point. Observe that the spherical suspension $\Sigma(\rrr P^2)$ of $\rrr P^2$ is an Alexandrov space, although not a manifold. Suspending the isometric $S^1$ action on $\rrr P^2$ yields an isometric $S^1$ action on $\Sigma(\rrr P^2)$ with fixed point set  a closed interval. 
 \end{example}

We now proceed with the proof of Theorem B. 
Let $G$ be a compact Lie group acting  effectively and by cohomogeneity one on a closed,  $3$-dimensional Alexandrov space $X^3$. Recall that  if the maximal dimension of an orbit is $k$ then $\dim(G)\leq k(k+1)/2$ (cf. \cite{MZ}). Since the dimension of the isometry group must be either $2$ or $3$ and the action is effective, $G$ must be $T^2$ or $SO(3)$.  To prove Theorem~B, we will show that the orbits of a $T^2$ or an $SO(3)$ action are 
the same as they are in the manifold case (cf. Proposition \ref{P:THMA} below), noting that
 in the  case where $G=SO(3)$, instead of obtaining only manifolds, 
  the group diagram given by $(SO(3), O(2), SO(3), SO(3))$ yields $\Sigma (\rrr P^2)$. Moreover, this is the only possibility up to equivariant homeomorphism since $N(O(2))/O(2)=\{e\}$. All other cases give us the same group diagrams as in the manifold case and it then follows that $X^3$ is one of the manifolds listed in the work of Mostert \cite{Mo} and Neumann \cite{N}.
  
  We further remark that one can apply theorem A to obtain the result by a straightforward analysis of 
  all possible group diagrams. We present here an alternative proof for the $T^2$ case that relies on the description of the space of directions at any point in $X^3$.

 
 \begin{proposition}\label{P:THMA} Let $X^3$ be a closed, $3$-dimensional Alexandrov space and suppose that $G$ acts isometrically on $X^3$ by cohomogeneity one.
 \begin{itemize}
 	\item[(1)]  If $G=T^2$, then $X^3$ is always a manifold.
	\\
  
  	\item[(2)] If $G=SO(3)$,  then there is one $X^3$ that is not a manifold, namely when the principal orbit is $SO(3)/O(2)$ and $X^3$ is $\Sigma(\rrr P^2)$. 
	\end{itemize}
 \end{proposition}

 \proof
 We recall that in dimension $3$ the space of directions at a point can be either  $S^2$ or $\rrr P^2$. Since $\rrr P^2$ is not the join of $S^1$ and a finite number of points, nor is it the join of $S^0$ and a $1$-dimensional set, it is clear that $\Sigma_x X^3$ can only be $\rrr P^2$ when $x$ is a fixed point of $G$. Since $T^2$ does not act 
 transitively on $\rrr P^2$, we need only consider the case where $\Sigma_x X^3=\rrr P^2$ when we consider the group $G=SO(3)$. In particular, this tells us in the case where $G=T^2$ that the space of directions at any point is an $S^2$, and thus $X^3$ itself is a manifold.

Now we consider the case where $G=SO(3)$. Since the action is of cohomogeneity one, the principal isotropy group contains $SO(2)$. Since the only infinite proper subgroups of $SO(3)$ are $SO(2)$ and $O(2)$, we see immediately that the only isotropy groups are $SO(3)$, $SO(2)$ and $O(2)$.
It follows  that the only case in which the space of directions is not $S^2$ is when the principal isotropy is $O(2)$ and the singular orbits are fixed by the $SO(3)$ action.
\qed


\section{Proof of Corollary B}
\label{Proof_Cor_B}
We now prove the corollary to Theorem B. We observe once again that it is sufficient to show that $X^n$ has the same orbit structure as a $T^{n-1}$ cohomogeneity one manifold of dimension $n$. It then follows  that a $T^{n-3}$ subgroup will act freely on $X^n$ (cf. \cite{P}) and thus $X^n$ is a $T^{n-3}$ bundle over $X^3$, a closed, $T^2$-cohomogeneity one Alexandrov space. Since the principal isotropy group is the trivial subgroup $\{e\}$, and $N(\{e\})=T^{n-1}$ is connected, it follows that 
the bundle must be trivial.

\begin{proposition} Let $T^{n-1}$ act by cohomogeneity one on $X^n$, a closed Alexandrov space. Then the only orbits possible are $T^{n-1}, T^{n-1}/\zzz_2$ and $T^{n-2}$.
\end{proposition}

\proof
A torus $T^k$ can only act effectively and 
transitively
 on a sphere when $k=1$ and the sphere is of the same dimension. The only finite subgroup of $T^k$ that can act effectively and 
 transitively on a sphere is $\zzz _2$ and the sphere is then $S^0$.  In particular, this tells us that the 
isotropy groups of the singular orbits can only be $T^{n-1}/\zzz_2$ or $T^{n-2}$. 
\qed


\section{Closed 4-dimensional Alexandrov spaces of cohomogeneity one}
\label{S:D4}

We first note that the possible groups that can act (almost) effectively by cohomogeneity one on $X^4$ 
are $SU(2)$, $SO(3)$, $T^3$, $SO(4)$ and $SO(3)\times S^1$. We will refer to a fiber bundle with fiber a cone over a homogeneous manifold $G/H$ as a \emph{$C(G/H)$-bundle}. We will denote a $C(G/H)$-bundle with base $B$ by $C(G/H) [B]$, and will denote a disc bundle over $Y$ by $D(Y)$.


\begin{thmc}Let $G$ act (almost) effectively and isometrically by cohomogeneity one on $X^4$, a closed Alexandrov space of dimension $4$. If $X^4$ is not a manifold, then it is given by one of the following cohomogeneity one group diagrams:
{\small
\begin{align*}
	&	(SU(2), \zzz_n, SU(2), SU(2))=\Sigma(SU(2)/\zzz_n),								\\
	&	(SU(2), \zzz_n, SU(2), \zzz_{2n})=C(SU(2)/\zzz_n)\cup_{SU(2)/\zzz_n}D(SU(2)/\zzz_{2n}),		
	\\
	&	(SU(2), \zzz_n, SU(2), S^1)= C(SU(2)/\zzz_n)\cup_{SU(2)/\zzz_n}D(S^2),					\\
	& (SU(2), <i>, SU(2), \{e^{j\theta}\}\cup \{ie^{j\theta}\})=C(SU(2)/<i>)\cup_{SU(2)/<i>}D(S^2),
	\\
	&	(SO(3), \Gamma, SO(3), SO(3))=\Sigma(SO(3)/\Gamma),								\\
	&	(SO(3), \zzz_n, SO(3), \zzz_{2n})=C(SO(3)/\zzz_n)\cup_{SO(3)/\zzz_n} D(SO(3)/\zzz_n),		\\
	& 	(SO(3), \zzz_n, SO(3), D_n)=C(SO(3)/\zzz_n)\cup_{SO(3)/\zzz_n} D(SO(3)/D_n)	,			\\
	& 	(SO(3), \zzz_n, SO(3), SO(2))=C(SO(3)/\zzz_n)\cup_{SO(3)/\zzz_n} D(S^2)	,				\\
	& 	(SO(3), D_n, SO(3), D_{2n})=C(SO(3)/D_n)\cup_{SO(3)/D_n} D(SO(3)/D_{2n}),				\\
	&	(SO(3), D_n, SO(3), O(2))=C(SO(3)/D_n)\cup_{SO(3)/D_n}  D(\rrr P^2),					\\
	&	(SO(3), T, SO(3), O)=C(SO(3)/T)\cup_{SO(3)/T}  D(SO(3)/O),							\\	
	&	(SO(3)\times S^1, O(2), SO(3), SO(3))=\Sigma(\rrr P^2)\times S^1,						\\
	&	(SO(3)\times S^1, O(2), SO(3), O(2)\times \zzz_2)=C(\rrr P^2)[S^1]\cup_{\rrr P^2\times S^1}D(\rrr P^2\times S^1),	\\
	&	(SO(3)\times S^1, O(2), SO(3), O(2)\times S^1)=C(\rrr P^2)[S^1]\cup_{\rrr P^2\times S^1} D(\rrr P^2),	 or		\\
	&	(SO(4), O(3), SO(4), SO(4))=\Sigma(\rrr P^3),
\end{align*}
}
where $\Gamma$ is one of $\zzz_n$, $D_n (n\geq 2)$, $T$, $O$ or $I$, the cyclic, dihedral, tetrahedral, octahedral or icosahedral subgroups of $SO(3)$, respectively and $Q=\{\pm 1, \pm i, \pm j, \pm k\}$.
\end{thmc}

\begin{remark} It will follow from the proof of Theorem~C that, when $X^4$ is not a cohomogeneity one $4$-manifold, it will be the suspension of a $3$-dimensional homogeneous spherical space form,  the union of a cone over a spherical space form and a disc bundle over a singular orbit, with certain restrictions,  or the union of a $C\rrr P^2$-bundle over $S^1$ and a disc bundle over a singular orbit, again with certain restrictions.
\end{remark}

\proof[Proof of Theorem~C] We proceed with a case by case analysis. Since we treated the case $G=T^3$ in the previous section, we will not repeat it here. Recall that the positively curved homogeneous spaces of dimension less than or equal to 3 are $S^3$, $S^3/\Gamma$, $S^2$, $\rrr P^2$, $S^1$ and $S^0$, where $\Gamma$ is one of 
$\zzz_n$, $n\geq 1$, $D_{m}^*$, $m\geq 2$, $T^*$, $I^*$, $O^*$ and $^*$ denotes the corresponding binary group; $T$, $I$, $O$ are, respectively, the tetrahedral, icosahedral and octahedral groups (cf. \cite{W}). By Theorem~A, $K_{\pm}/H$ are positively curved homogeneous spaces.
We need not consider the cases where $K_{\pm}/H$ are both spheres, as those were classified by Parker.
It remains to consider the cases where $K_{-}/H$ is a $3$-dimensional homogeneous spherical space form or $\rrr P^2$.

We note first that for the groups $SO(3)$ and $SU(2)$, the principal isotropy group is a finite group.
In this situation $K_{-}/H$ can be one of the $3$-dimensional spherical space forms and never $\rrr P^2$. When $G=SO(3)\times S^1$, the principal isotropy group must contain a circle and it follows directly from proposition~\ref{p:red} that it must be contained in the $SO(3)$ factor. Therefore $K_{-}/H$ is a positively curved 
 homogeneous space not homeomorphic to a sphere only when it is of dimension 
equal to $2$. In particular, $(SO(3)\times S^1)/(O(2)\times S^1)=\rrr P^2$. Finally when $G=SO(4)$, the principal isotropy group  is $3$-dimensional and thus $SO(3)$ or $O(3)$. Observe that $SO(4)/SO(3)=\rrr P^3$ is the only $3$-dimensional homogeneous space form that can occur and we never get an $\rrr P^2$.

We remark  that we follow Parker's convention in \cite{P}: any discrete group $H\subset SU(2)$ of even order
contains the central $\zzz_2$ and the only effective actions of $SU(2)$ have principal isotropy 
$\zzz_{2n+1}$. However, Parker views some actions as almost effective actions of 
$SU(2)$, rather than as effective actions of $SO(3)$, to unify the classification. In particular, 
he lists actions with principal isotropy 
$H=\zzz_n$ and $K_{\pm}\in \{\,\zzz_{2n},S^1,SU(2)\,\}$, for $n$ even or odd, as 
$SU(2)$ actions, regardless of the parity of $n$. He then lists all other almost effective actions of $SU(2)$ as $SO(3)$ actions.

We further remark that there was an omission in Parker~\cite{P}, which was pointed out by Hoelscher in his thesis~\cite{H}, where he gives a classification of low-dimensional simply-connected cohomogeneity one manifolds.  Namely, the group diagram 
\bdm 
(SU(2), <i>, \{e^{i\theta}\}, \{e^{j\theta}\cup \{i e^{j\theta}\}),
\edm 
corresponding to the $SO(3)$ action on $\ccc P^2$ viewed as a sub-action of the homogeneous 
action of $SU(3)$, is missing from the simply-connected classification. This means that the group diagrams 
\bdm
(SU(2), <i>,\{e^{j\theta}\}\cup \{ie^{j\theta}\}), \{e^{j\theta}\}\cup \{ ie^{j\theta}\}),
\edm
\bdm 
(SU(2), <i>, Q, \{e^{j\theta}\}\cup \{ ie^{j\theta}\})
\edm
and
\bdm 
(SU(2), <i>, <e^{\frac{i\pi\theta}{4}}>, \{e^{j\theta}\}\cup \{i e^{j\theta}\}),
\edm
where $Q=\{ \pm 1, \pm i, \pm j, \pm k\}$, 
 are therefore also missing from the manifold classification. All other possibilities arising from the group diagram corresponding to the simply-connected case already form a part of the classification given by Parker. Thus we will have one more possibility for $K_{-}$
when we consider the group $SU(2)$, namely, in the case where $H=<i>\cong \zzz_4$, $K_-=\{e^{j\theta}\}\cup \{ie^{j\theta}\}$.

We now proceed with the case by case analysis.
\\


\noindent \textbf{Case: ${\bf G=SU(2)}$.} 
As mentioned above, the principal orbit is $3$-di\-men\-sional and thus the principal isotropy group $H$ can be trivial or $\zzz_n$.  The only possibilities for  $K_{\pm}/H$ that are not spheres are the homogeneous spherical space forms $SU(2)/\zzz_n$.

Clearly, $\Sigma(SU(2)/\zzz_n)$, the suspension of $SU(2)/\zzz_n$, admits a cohomogeneity one action. The group diagram is given by $(SU(2), \zzz_n, SU(2), SU(2))$.
We now consider the cases when $K_{-}/H=SU(2)/\zzz_n$ and $K_{+}/H$ is a sphere.
The only other possibilities for the subgroups $K_{+}$ are $\zzz_{2n}$, $S^1$ and a $\zzz_2$ extension of $S^1$ (the last case only occurs in the omitted case pointed out by Hoelscher, where $H=<i>\cong \zzz_4$, see the remarks above), and thus we obtain
an Alexandrov space that decomposes as the union of a cone over $SU(2)/\zzz_n$ and a disc bundle over one of $SU(2)/\zzz_{2n}$, $SU(2)/S^1=S^2$ or $SU(2)/(S^1\times \zzz_2)=S^2$. The group diagrams are given by $(SU(2), \zzz_n, SU(2), \zzz_{2n})$, $(SU(2), \zzz_n, SU(2), S^1)$ and $(SU(2), <i>, SU(2), \{e^{j\theta}\}\cup  \{ie^{j\theta}\})$, respectively.
\\

\noindent \textbf{Case: ${\bf G=SO(3)}$.}
The principal isotropy group $H$ can be trivial, $\zzz_n$, $D_{m}, m\geq 2$, $T$, $I$ or $O$. We obtain the spherical suspension of $\rrr P^3$ or of the homogeneous spherical space forms $SO(3)/\Gamma$, where $\Gamma$ is one of $\zzz_n$, $D_{m}, m\geq 2$, $T$, $I$ or $O$. The group diagrams are given by $(SO(3), \Gamma, SO(3), SO(3))$. The remaining possibilities 
are when $K_{-}/H$ can be one of the previous homogeneous spaces and $K_{+}/H$ is a sphere.
In this case, the principal isotropy group $H$ is either the trivial group, $\zzz_n, D_{2m}$ or $T$.
We can decompose the Alexandrov space as follows, depending on the possibilities for $H$. 

\begin{enumerate}
	\item ${\bf H=\zzz_n}$. Here $X^4$ will be the union of the cone over $SO(3)/\zzz_n$ and a disc bundle over one of $SO(3)/\zzz_{2n}$, $SO(3)/D_n$  or 
	$SO(3)/SO(2)=S^2$. The group diagrams are given by 
	 $(SO(3), \zzz_n, SO(3), \zzz_{2n})$, $(SO(3), \zzz_n, SO(3), D_n)$ and 
	 $(SO(3), \zzz_n, SO(3), SO(2))$, respectively.

	\item ${\bf H=D_n}$. In this case $X^4$ will be the union of the cone over $SO(3)/D_n$ and a disc bundle over $SO(3)/D_{2n}$ 
	or a disc bundle over $SO(3)/O(2)=\rrr P^2$. The group diagrams are given by $(SO(3), D_n, SO(3), D_{2n})$ and $(SO(3), D_n, SO(3), O(2))$, respectively.

	\item ${\bf H=T}$. Here $X^4$ will be the union of the cone over $SO(3)/T$ and a disc bundle over $SO(3)/O$. Its group diagram is given by $(SO(3), T, SO(3), O)$.
\end{enumerate}

\noindent \textbf{Case:  ${\bf G=SO(3)\times S^1}$.} For this group, the principal isotropy subgroup is one-dimensional and hence contains a circle subgroup, $T^1$. By proposition~\ref{p:red}, it is clear that if $\pi_2(T^1)=S^1$, that is, the projection of  this circle subgroup onto the second factor is the second factor, then the action is reducible to an $SO(3)$ cohomogeneity one action, which is treated above. Thus, we need only consider the case where the principal isotropy must be entirely contained in the first factor and there are only two possible principal isotropy subgroups: $SO(2)$ and $O(2)$.
It is then immediate that the only principal isotropy group for which we obtain a non-spherical fiber over a singular orbit  is $O(2)$.

When both singular orbits have isotropy group $SO(3)$, $X^4$ will decompose as the union of $C\rrr P^2$-bundles over the singular orbits $(SO(3)\times S^1)/SO(3)=S^1$. By Corollary  \ref{c:interval}, we see that this is a $\Sigma(\rrr P^2)$-bundle over $S^1$.
Further, $N(H)/H=S^1$ is connected and so the bundle is trivial.
The group diagram is given by $(SO(3)\times S^1, O(2), SO(3), SO(3))$.

The only other case to consider is when $K_{+}/H$ is a sphere.
In this case, $X^4$ will decompose as the union of a $C\rrr P^2$-bundle over the singular orbit $(SO(3)\times S^1)/SO(3)=S^1$, and a  disc bundle over one of $(SO(3)\times S^1)/(O(2)\times \zzz_2)=\rrr P^2\times S^1$ or $(SO(3)\times S^1)/(O(2)\times S^1)=\rrr P^2$.
The group diagrams are given by $(SO(3)\times S^1, O(2), SO(3), O(2)\times \zzz_2)$ and $(SO(3)\times S^1, O(2), SO(3), O(2)\times S^1)$, respectively.
\\

\noindent \textbf{Case: ${\bf G=SO(4)}$.} In this case the only possibility is that $X^4$ is $\Sigma(\rrr P^3)$. The group diagram is given by 
$(SO(4), O(3), SO(4), SO(4))$.

\qed

\begin{remark} It is worth noting that in dimension $3$ the suspension of $\rrr P^2$ is the only Alexandrov space of cohomogeneity one that is not a manifold. In dimension $4$ we 
have suspensions of spherical space forms,  $S^2$-bundles over the suspension of $\rrr P^2$, and the union of disc bundles with cones over spherical space forms or bundles over a circle with fiber  $C(\rrr P^2)$. However we never obtain the union of two $C(K_{\pm}/H)$-bundles over different singular orbits.
This phenomenon will occur in higher dimensions, however. Consider, for example, the group diagram
$$(SO(3)\times SO(3), O(2)\times O(2), SO(3)\times O(2), O(2)\times SO(3)),$$
 corresponding to a cohomogeneity one action on an Alexandrov space of dimension $5$. This is the union of two $C(\rrr P^2)$-bundles over distinct $\rrr P^2$. \end{remark}


\end{document}